\documentclass[draft]{agujournal2019}
\usepackage{url} 
\usepackage[inline]{trackchanges} 
\usepackage{soul}

\usepackage{amsmath,amssymb,bm}
 
\draftfalse

\journalname{Journal of Advances in Modeling Earth Systems (JAMES)}

\begin{document}

\title{Data Assimilation in Chaotic Systems Using Deep Reinforcement Learning}

%
%

\authors{Mohamad Abed El Rahman Hammoud\affil{1}, Naila Raboudi\affil{1}, Edriss S. Titi\affil{2,3}, Omar Knio\affil{1} and Ibrahim Hoteit\affil{1}}

\affiliation{1}{King Abdullah University of Science and Technology, Thuwal 23955, Saudi Arabia}
\affiliation{2}{Department of Applied Mathematics and Theoretical Physics, University of Cambridge, Cambridge CB3 0WA, UK}
\affiliation{3}{Department of Mathematics, Texas A \& M University, College Station, TX 77843, USA}

\correspondingauthor{Ibrahim Hoteit}{ibrahim.hoteit@kaust.edu.sa}

\begin{keypoints}
\item Deep reinforcement learning (RL) is introduced for data assimilation 
\item RL generalizes to new situations unseen during training through actively learning from the data and system dynamics
\item The RL agent allows for nonlinear state-adaptive correction of the forecast using the observations
\item The performance of the proposed RL algorithm surpasses that of the ensemble Kalman filter (EnKF) with the Lorenz '63
\end{keypoints}

\begin{abstract}
Data assimilation (DA) plays a pivotal role in diverse applications, ranging from climate predictions and weather forecasts to trajectory planning for autonomous vehicles. 
A prime example is the widely used ensemble Kalman filter (EnKF), which relies on linear updates to minimize variance among the ensemble of forecast states. 
Recent advancements have seen the emergence of deep learning approaches in this domain, primarily within a supervised learning framework. 
However, the adaptability of such models to untrained scenarios remains a challenge.
In this study, we introduce a novel DA strategy that utilizes reinforcement learning (RL) to apply state corrections using full or partial observations of the state variables. 
Our investigation focuses on demonstrating this approach to the chaotic Lorenz '63 system, where the agent's objective is to minimize the root-mean-squared error between the observations and corresponding forecast states. 
Consequently, the agent develops a correction strategy, enhancing model forecasts based on available system state observations.
Our strategy employs a stochastic action policy, enabling a Monte Carlo-based DA framework that relies on randomly sampling the policy to generate an ensemble of assimilated realizations. 
Results demonstrate that the developed RL algorithm performs favorably when compared to the EnKF. 
Additionally, we illustrate the agent's capability to assimilate non-Gaussian data, addressing a significant limitation of the EnKF.
\end{abstract}

\section*{Plain Language Summary}

\noindent Reliable forecasts of the state of chaotic systems, such as environmental flows, require combining observational data and dynamical model outputs through a process called data assimilation. 
The ensemble Kalman filter (EnKF) is the most commonly adopted algorithm for this task, however, is subject to some limitations when applied to nonlinear/non-Gaussian systems. 
Recently, there has been interest in using deep learning (DL), particularly within a supervised learning setup, for DA.
However, making DL models work well in new situations that differ from those experienced during training is challenging.
In this work, we propose a new DA approach that leverages reinforcement learning (RL). 
RL helps the system make corrections to its predictions based on observed data, even if the model hasn't been trained for those specific scenarios. 
Compared to the state of the art DA algorithms, RL offers a novel framework for nonlinear corrections of the forecast using the incoming observations.
Numerical results show that the proposed RL algorithm outperforms the EnKF and demonstrates the RL agent's ability at addressing some shortcomings of the EnKF.

\section{Introduction}

Assimilating observational data is essential for improving predictability and understanding complex dynamics in chaotic and dynamic physical systems. 
Chaotic dynamical systems, such as those describing climate and weather, involve inherent imperfections and extreme sensitivity to initial conditions, whereas the observational data available for such systems often carry significant uncertainties \cite{Eckmann1985}.
To address the associated challenges, data assimilation (DA) combines real-world observations with numerical model outputs, continually refining model predictions by aligning them with newly acquired observations to enhance the accuracy and reliability of the predictions \cite{Kalnay2004}.
DA techniques are broadly categorized as variational and filtering methods \cite{le1986variational,Ghil91, Lorenc2003,hoteitchapter2018}. 
The ensemble Kalman filter (EnKF) represents one of the most popular filtering DA techniques, especially in the context of large-scale nonlinear systems \cite{Evensen2003}. 
Operating within a Bayesian probabilistic framework, the EnKF squentially splits the filtering (state estimation) process into cycles that alternate between forecast steps, driven by the system's dynamical model, and analysis steps, which updates the forecast with incoming data \cite{Evensen2003}.
This approach enables Gaussian-based Monte Carlo (MC) approximations of both state forecast and analysis distributions through an ensemble of state samples \cite{Hoteit2008}. 

EnKF schemes are considered as the gold standard when assimilating uncertain observations of the system states across diverse fields due to their robustness, capacity to handle complex and high-dimensional systems, and computational efficiency \cite{Houtekamer98_EnKF_DA}. 
However, their applicability is not without constraints, particularly when the underlying assumptions are compromised. 
In particular, challenges may arise from the EnKF's inherent linear assumption, and the necessity for maintaining a Gaussian distribution within the ensemble, both of which become challenging in the presence of strong nonlinearities \cite{kalnay_2002,Hoteit2008}. 
Additionally, whereas the Gaussian assumption for both model and observational noise offers convenience, it may not universally hold in real-world scenarios, thereby limiting EnKF's performance, especially when errors deviate significantly from Gaussian patterns. 
In such cases, it is necessary to explore alternative approaches that are better suited for these scenarios; e.g.~\citeA{PJVL2009}.

Reinforcement Learning (RL) is a paradigm of artificial intelligence that deals with how an agent can learn to make decisions through interactions with an environment, namely to achieve a specific objective \cite{Recht2019}.
It is inspired by behavioral psychology and focuses on learning how to take actions in an environment to maximize some notion of cumulative reward.
Within the RL framework, an agent engages in trial-and-error exploration, testing various actions and observing their outcomes \cite{Mnih2015}.
The agent's goal is to formulate an optimal strategy, often referred to as a policy, that guides its actions to maximize the cumulative reward over a time horizon.
It is noteworthy to point out that the RL framework is inherently different from the supervised learning approaches because the latter require a pre-computed reference database for training, which in this context consists in minimizing a global objective function \cite{Glorot2010, Karniadakis2021}.
RL finds extensive applications in domains necessitating dynamic control and decision-making capabilities, encompassing fields such as robotics \cite{Kober2013}, gaming \cite{mnih2013, Vinyals2019}, autonomous navigation \cite{Sallab_2017}, fluid dynamics \cite{Novati2021, Koumoutsakos2022}, and beyond.

In this work, we introduce a novel DA formalism utilizing RL as a strategy to actively update a nonlinear forecast correction scheme with the incoming data.
The RL agent learns through interactions with the environment, adapting to its changes, and actively applies nonlinear corrections to handle complex processes.
Numerical experiments were conducted with the Lorenz '63 chaotic system \cite{Lorenz63}, and the RL agent's performance was benchmarked against the EnKF algorithm using a large cardinality ensemble under various experimental conditions.
These include tracking a reference solution and assimilating normally-distributed noisy observations at various noise levels and observation frequencies. 
Furthermore, we investigate the performance of the RL agent at assimilating observations with different noise distribution models, namely uniform, log-normal and Gaussian noise. 
We further explore the RL agent's effectiveness at assimiliating partial state observations.

The remaining of the manuscript is organized as follows. 
Section \ref{intro} introduces the RL-DA framework.
The RL methodology for DA is then described in Section \ref{methods}, where a comprehensive overview of the RL framework is first introduced, accompanied by a description of the Lorenz '63 system and the EnKF algorithm. 
Sections \ref{results1} and \ref{results2} present our numerical results. 
Finally, Section \ref{conclusion} summarizes the main conclusions of this study.

\begin{figure}[!htbp]
    \centering
    \includegraphics[width=0.8\linewidth]{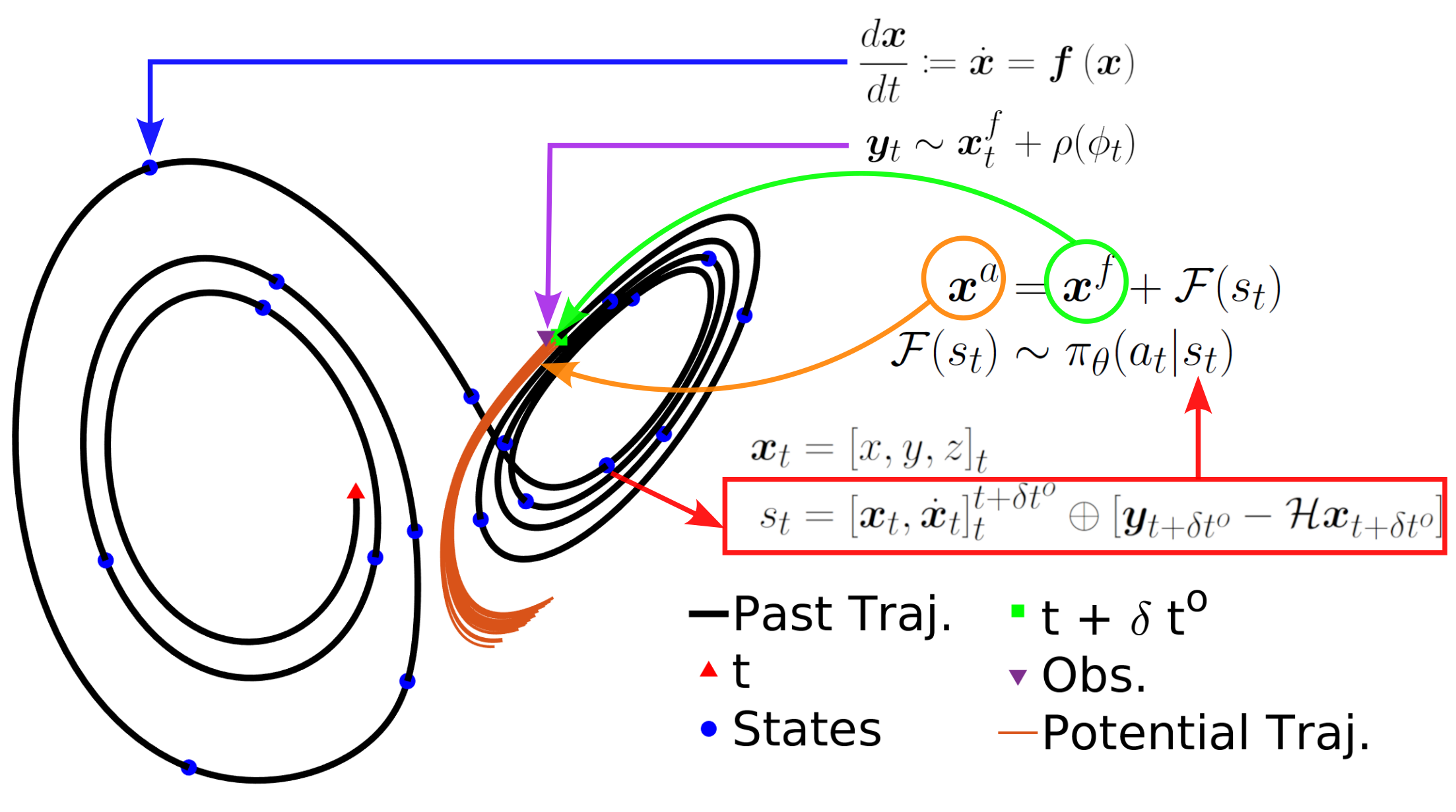}
    \caption{Schematic of the proposed reinforcement learning-based data assimilation framework using the Lorenz '63 as the main example. The plot illustrates the Lorenz '63 solution trajectory (black curve) with an arbitrary assimilation window start time $t$ (red triangle) and corresponding end time $t+\delta t^o$ (green square) when a new observation is available and assimilated. The three dimensional state variables ($\bm{x}$) of the model are shown at every model time step $\delta t$ (blue circles). At the last time step, the noisy observational data point ($\bm{y}$) is shown (inverted purple triangle) alongside the different evolution trajectories (orange curves) following several corrections ($\mathcal{F}(s_t)$) sampled from the policy function $\pi_{\theta}(a_t | s_t)$. 
    The policy $\pi_{\theta}(a_t | s_t)$ considers as input state vector the extended state vector composed of the concatenation of the forecast state variables ($\bm{x}$) and their time derivatives ($\dot{\bm{x}}$) at each time step $\delta t$ between $t$ and $t + \delta t^o$ alongside the innovation term, defined as the difference between the observation and its correspondent forecast. The concatenation operation is denoted by $\oplus$, and for the sake of conciseness, concatenation of $\bm{x}$ and $\dot{\bm{x}}$ at each $\delta t$ is represented by the sub- and super-scripts of $\left[ \bm{x},\dot{\bm{x}} \right]$.
    Since a stochastic policy is considered in the DA framework, an ensemble of $\mathcal{F}(s_t)$ correction terms are sampled from $\pi_{\theta}(a_t | s_t)$ when a noisy observation is available.
    Note that the state variables might not be fully observed, hence $\mathcal{H}$ projects the forecast onto the observation space. Moreover, the observation $\bm{y}$ is considered to be a noisy estimate of the forecast with no restriction on the distribution of the additive noise.}
    \label{fig:graphicalAbstract}
\end{figure}

\section{Reinforcement Learning For Data Assimilation}
\label{intro}

In RL, agents make sequential decisions to achieve specific goals, with the focus on maximizing cumulative rewards over time \cite{SuttonBarto, bertsekas2019}. 
This aligns with decision-making scenarios where actions have consequences, and objectives must be met. 
RL is particularly relevant to control systems \cite{Titi2014, Titi2018}, where agents learn control policies to influence the behavior of systems \cite{Silver2014}. 
The key concept in RL is the trade-off between exploration, where the agent experiments with new actions, and exploitation, where the agent chooses known actions with high rewards, mirroring real-world decision-making challenges \cite{Sallab_2017}. 
RL agents learn from feedback, adapt to changing environments, and generalize knowledge to make decisions in new situations. 

DA is an essential process used in scientific fields such as meteorology, oceanography, and environmental modeling to guide the state of complex systems with incoming observations \cite{Ghil91, Hoteit2008}. 
It involves merging observational data with numerical models to enhance predictions once observational information is available \cite{kalnay_2002}. 
This process continuously drives the computed system state to align with observations, thereby ensuring accurate and robust state estimates.
DA accounts for model and observational uncertainties, offering more reliable predictions for chaotic systems, making it indispensable for tasks such as weather forecasting \cite{Rabier2005} and climate modeling \cite{Pedatella2014}.
Hence, adopting an RL framework for DA is a natural progression in the domain, enabling for a nonlinear correction scheme that is also free from restrictive assumptions on the statistics of the observations and model. 

In RL, an agent exists in an environment that is described by a set of dynamical rules characterizing its evolution, for example, a system of differential equations \cite{SuttonBarto}. 
The agent's responsibility is to make decisions affecting its environment in a way that it maximizes the cumulative reward, or achieves a particular goal. 
The ultimate outcome of the RL's training procedure is an agent policy $\pi_{\theta}(a_t | s_t)$, a mapping from the observation space to the action space, which is evaluated to actively control the behavior of the agent at state $s_t$ in a dynamical system.
The policy function is generally characterized by a neural network parameterized with $\theta$.
Policy functions can be categorized as either deterministic or stochastic; in a deterministic policy, the action with the highest probability is chosen, whereas a stochastic policy relies on random sampling to select an action.
In the present framework, a stochastic policy was adopted from which the DA correction term was sampled, where actions are sampled from a Gaussian policy \cite{Schulman2017}.
Hence, after training, a policy function is obtained that could be used to sample potential correction terms from a distribution that adapts to the agent's state, and allowing to generate an ensemble of states via MC sampling.
In contrast with most efforts put for developing efficient DA schemes; eg.\ \cite{LERMUSIAUX2007172, Farchi2021, Buizza2022}, the RL machinery relies on a nonlinear neural network to provide a correction without being restricted to a pre-computed dataset for supervising its training. 
Furthermore, the RL agent does not require any assumption on the noise distribution of the observational errors nor restrictive assumptions on the model. 

In this study, the chaotic Lorenz '63 system of differential equations was considered to examine the performance of RL at DA for a chaotic dynamical system \cite{Lorenz63}. 
The system describes the solution of a three-dimensional state vector, $\bm{x} = \left[x, y, z\right]$; it is characterized by a chaotic attractor, where the solution is sensitive to initial conditions and experiences a nonperiodic behavior \cite{Eckmann1985, Bakarji2023}.
In this setting, the agent receives information, in the form of an extended state vector describing the system, denoted by states, that includes the forecast states and their derivatives $\bm{x}^f$ and $\dot{\bm{x}}^f$, respectively, at each model time step $\delta t$ starting from the time $t$ of the previous observation till the next observational time step $t +\delta t^o$, and the innovation term $\bm{y} - \mathcal{H}\bm{x}^f$.
Here, $\mathcal{H}$ represents the observation operator that projects the model forecast $\bm{x}$ onto the observation space and $\bm{y}$ denotes a noisy observation of the system state. 

The agent interacts with the environment to change its course of evolution and adapts to these changes to maximize the cumulative reward, as later defined, gathered over some period of time \cite{Silver2014}. 
This interaction is formulated mathematically as:

\begin{equation}
    \bm{x}^a = \bm{x}^f + \mathcal{F} \left( s_t \right),
\end{equation}
where the corrected state vector, denoted by a superscript $a$ for analysis, $\bm{x}^a$ is the sum of the model forecast, $\bm{x}^f$, and the correction term $\mathcal{F} \left( s_t \right)$, which is sampled from $\pi_{\theta} (a_t |s_t)$.
Note that this form of the update is similar to that of the Kalman Filter and the EnKF algorithm, however, the latter rely on a linear update term \cite{Kalman1960}. 
In the current configuration, the RL agent is not provided with statistical information regarding the noisy observations. 
Instead, it employs an MC strategy, using an RL agent that employs random stochastic policy sampling. 
This approach generates an ensemble of assimilated solutions, which are subsequently averaged to produce an improved estimate of the system state, denoted by RL-50 in the following sections.

The training cycle is defined by specifying the reward function \cite{lillicrap2015}. 
We test out several reward functions in our preliminary investigation, where the agent's performance was evaluated using the mutual information, negative of the root-mean-squared error (RMSE) and $\text{RMSE}^{-1}$. 
While these reward functions are mathematically similar \cite{Seidler1971, Guo2005}, the associated training stability is different.
Accordingly, the agent was trained to maximize the negative of the RMSE, which strikes a satisfactory balance between interpretability, computational cost and agent's performance.

\section{Methods}
\label{methods}

\subsection{Reinforcement Learning}

The framework for RL involves training an agent through several interactions with the environment, in the present context, the dynamical system. 
Training an RL agent requires a large number of interactions with the environment and consequently a large unavoidable computational load often several orders of magnitude greater than solving the underlying differential equations. 
However, the field of RL has become more accessible in recent times, thanks to open-source libraries like \texttt{smarties} \cite{novati2019a} and \texttt{stable baselines3} \cite{raffin2021stable}, among others. 
In this work, we leverage the capabilities of \texttt{stable baselines3}, a high-performance RL software designed to exploit parallel computing, distributing the training process across multiple computational nodes. 
In the present configuration, each node simulates a distinct trajectory of the Lorenz '63 system, providing a large set of agent-environment interactions that are used to train the agent.
In this parallelized setup, each computational node accumulates experiences by independently interacting with various instances of the environment. 
These experiences are then structured into episodes defined as:

\begin{equation}
    \tau = \{ s_t, r_t, a_t, s_{t+1}\}_{0:T},
\end{equation}
where $\tau$ is the ordered set of interactions across a time horizon, t represents the time at which the environment is at state $s_t$, $a_t$ is the action the agent takes at that time, $r_t$ is the reward the agent receives from performing action $a_t$ and $s_{t+1}$ is the subsequent state. 

The RL agent's training objective is to maximize the expected cumulative discounted reward function, defined as:

\begin{equation}
    R_t = \sum_{t=0}^{T} \gamma^{t} r_t,
\end{equation}
where $\gamma \in [0, 1)$ is the discount factor. 
In our specific setting, a smaller value of $\gamma$ proves advantageous, given the random noise sampling. 
This choice of reducing the emphasis on distant future rewards results in a more stable agent performance.

The policy function $\pi_{\theta}$ is a mapping between the agent's state and the action space, which can be structured either as a set of discrete actions or as a probability distribution function for continuous actions. 
As previously mentioned, policy functions are either deterministic, the action to most likely result in the highest reward is chosen, or stochastic, where actions are randomly sampled from a distribution that is typically approximated by a surrogate model.
Here, the policy $\pi_{\theta}$ is represented as a densely connected multi-layer perceptron \cite{chen1995universal} parameterized by $\theta$.
Furthermore, actions assume continuous values, leading $\pi_{\theta}$ to output a probability distribution over possible actions. 
Hence, the agent's actions can be either sampled from this distribution, allowing the agent to explore the environment and seek potentially rewarding outcomes, otherwise, the action with the highest probability can be chosen.

\subsection{Proximal Policy Optimization}

In the present framework, we adopt the Proximal Policy Optimization (PPO) algorithm \cite{Schulman2017} and briefly describe it here for completeness. 
PPO trains an agent using two key components, each parameterized by distinct neural networks: an actor network that takes the environment's state as input and produces the corresponding action, and a critic network that also takes the environment's state as input and predicts the discounted reward \cite{Mnih2016}. 
In our study, both the actor and critic networks are represented by multi-layer perceptrons, each composed of two hidden layers, each containing 128 neurons.

The essence of the PPO algorithm revolves around optimizing the actor network to maximize the cumulative reward obtained by the agent, and the critic network to minimize the mean squared error between the predicted and actual expected cumulative rewards, starting from a given state. 
This optimization can be mathematically expressed through two distinct loss functions.
The actor network is optimized by maximizing the actor's objective function:

\begin{equation}
    J_{actor} = \mathbb{E} \left[ \text{min}\left(q_t \left(\theta\right) \hat{A}_t, \text{clip} \left( q_t\left(\theta \right), 1 - \epsilon, 1+\epsilon \right) \hat{A}_t \right) \right],
\end{equation}
where $q_t(\theta) = \pi_{\theta}(a_t|s_t) / \pi_{old}(a_t|s_t)$ is the ratio of the probability of adopting an action $a_t$ at state $s_t$ using $\pi_{\theta}$ to that of the previous policy $\pi_{old}$.
Furthermore, the present setting relies on policy clipping with an $\epsilon = 0.2$ \cite{Schulman2017}, where $q_t(\theta) \in \left[ 1-\epsilon , 1+\epsilon \right]$. 
This policy clipping mechanism helps maintain policy stability during parameter updates, stabilizing the training process.
On the other hand, the critic loss is given as:

\begin{equation}
    L_{critic} = \mathbb{E} \left[ \hat{A}^2  \right],
\end{equation}
where, $\mathbb{E}$ is the expectation operator and $\hat{A}$ is the advantage \cite{Mnih2016}, which quantifies how favorable the observed outcome of selecting a particular action is compared to the estimated discounted reward of the current state. 
The advantage is described as:

\begin{equation}
    \hat{A} = V_{target} - V_{\theta, old},
\end{equation}
where, $V_{target} = \sum_{i=0}^{T-1} r_{i} \gamma^{i} + \gamma^T V_{\theta, old}(s_{T})$ is the discounted reward computed using the agent's interactions with the environment and $V_{\theta, old}$ is the discounted reward predicted by the critic network.

\subsection{Lorenz '63}

The Lorenz '63 is a set of three deterministic ordinary nonlinear differential equations developed to simulate simplified atmospheric convection \cite{Lorenz63}.
This system is renowned for its manifestation of chaotic behavior, where even minuscule perturbations in initial conditions lead to substantially divergent solution trajectories over time \cite{Eckmann1985}.
The Lorenz equations have been extensively studied in chaos theory and nonlinear dynamics, and have been the fundamental benchmark to develop new data assimilation techniques \cite{Titi2001, Hayden2011}. 
The Lorenz '63 equations are given by:

\begin{align}
    \dot{x} &= \sigma(y - x),\\
    \dot{y} &= x(\rho - z) - y,\\
    \dot{z} &= xy - \beta z,
\end{align}
where, $\sigma$, $\rho$ and $\beta$ are typically positive constants. 
This system is known to exhibit a chaotic attractor for $\sigma = 10$, $\rho=28$ and $\beta=8/3$.
In this study, the system of equations were solved using an $2^{nd}$ order Runge-Kutta scheme with a time step $\delta t = 0.001$, which offers a suitable balance between solution accuracy and computational time for the application at hand.

\subsection{Data assimilation using Reinforcement Learning}

The present study explores a novel data assimilation framework that leverages RL to assimilate noisy observations of the system states and improve  estimates of the system states. 
In this investigation, the environment is represented by the chaotic Lorenz '63 system \cite{Lorenz63}. 
The RL agent receives noisy information about the system's state variables, and its policy, $\pi_{\theta}(a_t | s_t)$ that is contingent upon the environment's state $s_t$ takes an action according to the preassigned strategy.
The state upon which the agent's policy is evaluated consists of the extended vector composed by the concatenation $\left[\bm{x}^f, \dot{\bm{x}}\right]_{t}^f \oplus \left[\bm{x}^f, \dot{\bm{x}}\right]_{t + \delta t}^f \oplus ... \left[\bm{x}^f, \dot{\bm{x}}\right]_{t + \delta t^o}^f \oplus \left[ \bm{y}_{t+\delta t^o} - \mathcal{H}\bm{x}^f_{t+\delta t^o} \right]$. 
Notably, this selection preserves the Markovian assumption inherent in the EnKF, as $\mathcal{F} \left( \bm{x}|_{t}^{t+\delta t^o} \right) = \mathcal{F} \left( \bm{x}(t + \delta t^o) \right)$. 
However, including forecast information from previous steps significantly enhances training stability, even though it comes at the cost of a higher dimensional input. 
This gives rise to the question of how long back-in-time should forecast states be considered.

In this context, we introduce an RL agent responsible for correcting model forecasts of the dynamical system states using the update equation:

\begin{equation}
    \bm{x}_{t+\delta t^o}^a = \bm{x}_{t+\delta t^o}^f + \mathcal{F}_{\theta}\left(\bm{x}|_{t}^{t+\delta t^o}, \dot{\bm{x}}|_{t}^{t+\delta t^o}, \bm{y}_{t+\delta t^o} - \mathcal{H}\bm{x}_{t+\delta t^o}^f\right),
\end{equation}
where, $\mathcal{F}_{\theta}$ represents the RL agent's policy, parameterized by $\theta$. 
The policy takes as input the state vector $\bm{x}$ and the first-order derivatives $\dot{\bm{x}}$ at all time steps from $t$ to $t+\delta t^o$ at $\delta t$ increments, as well as the innovation term $\bm{y} - \mathcal{H}\bm{x}^f$. 
Since a stochastic policy function is considered, the study examines the performance of a single RL agent by taking maximum probability action, and the performance of an ensemble of agents by randomly sampling the policy function for actions.

\subsection{Training the DA agent}

The present experimental setup encompasses various hyper-parameters that require tuning to achieve satisfactory performance. 
The parameters subjected to tuning include the learning rate, $\gamma$, number of assimilation steps per episode ($n_{a, train}$), total number of episodes, value function coefficient ($v_f$), gradient clipping coefficient. 
Experiences have shown that the performance of a stable agent is most sensitive to $\gamma$, $v_f$ and gradient clipping. 

The process of hyper-parameter optimization commenced with a Latin hypercube sampling strategy to establish a baseline assessment of the acceptable range of values for these parameters. 
Subsequently, the training process is repeated using a new set of hyperparameters selected from within a finer-scale parameter space. 
For all experiments conducted, we employed the ADAM stochastic optimization algorithm \cite{Adam2017} to optimize the loss function for the parameters of the actor and critic networks.
The parameters utilized for training the agents, which underpin the results presented in this study, are detailed in Supplementary Table 1.

The RL agent's training objective centered on maximizing the cumulative rewards accrued over a specific time horizon. 
At each assimilation step, the reward was calculated as the negative RMSE between the observation and the forecast generated by the preceding action. 
This choice was made because minimizing the RMSE is equivalent to maximizing the mutual information between the compared quantities and because the RMSE is ultimately the measure that is used to evaluate the performance of the agent.
More specifically, since the experiments in this study feature a well-defined reference solution, we report the RMSE of both the RL and EnKF solutions with respect to the noise-free reference solution. 
The RMSE hence provides quantitative estimates that help examine the assimilated solution in terms of forecast and analysis.

\subsection{Ensemble Kalman Filter}

The EnKF algorithm is commonly employed to estimate a discrete-time state process, denoted as ${\bf{x}} = {\{{\bf{x}}_n\}}_{n \in \mathbb{N}}$, based on observations from a corresponding process ${\bf{y}} = {\{{\bf{y}}_n\}}_{n \in \mathbb{N}}$.
These processes are conventionally connected through a state-space system described as follows:

\begin{equation} 
\label{sys0} 
\left\{ 
\begin{array}{ccc} 
	\mathbf{x}_{t} &=& \mathcal{M} (\mathbf{x}_{t-1})  + \mathbf{u}_{t} \\
	\mathbf{y}_t &=& \mathcal{H} ({\bf{x}}_t) + \mathbf{v}_{t},  
\end{array}
\right. , 
\end{equation}
where $\mathcal{M}$ represents the nonlinear dynamical model, that advances the system state from time $t$ to $t+\delta t$, and $\mathcal{H}_t$ the observation operator that projects ${\bf{x}}_{t}$ from the state space onto the observation space.
Here, we make the simplifying assumption that $\mathcal{H}$ is linear, although EnKF algorithms can readily accommodate cases with nonlinear $\mathcal{H}$.
The noise terms, $\mathbf{u} = {\{\mathbf{u}_{t}\}}_{t \in \mathbb{N}}$ and $\mathbf{v} = {\{\mathbf{v}_{t}\}}_{t \in \mathbb{N}}$ are respectively the model and observation process noises.
The EnKF algorithm assumes $\mathbf{u}_t$ and $\mathbf{v}_t$ to follow Gaussian distributions with zero means and covariances ${\bf{Q}}_{t}$ and ${\bf{R}}_t$, respectively.
Furthermore, $\mathbf{u}$ and $\mathbf{v}$ are assumed to be independent, jointly independent and independent of the initial state ${\bf x}_0$.

The filtering problem involves estimating the state, ${\bf{x}}_t$, based on observations up to time $t$. 
EnKF algorithms are primarily designed to provide a MC approximation of the system state distribution using an ensemble of system state realizations.
From this ensemble, empirical estimates of the posterior mean state and associated error covariances are derived, typically in the form of sample means and covariances.
The process starts with an analysis ensemble of size $N_{ens}$ denoted as ${\{ { {\bf{x}}}_{t}^{a,i} \}}_{i=1}^{N_{ens}}$ available at time $t$. 
Subsequently, the forecast ensemble at the next time step $t+\delta t$ is computed by advancing each member ${{\bf{x}}}_{t-1}^{a,i}$ forward in time using the dynamical model, described as:

\begin{eqnarray}
{{{\bf{x}}}}_{t+\delta t}^{f,i} = \mathcal{M} ({{\bf{x}}}_{t}^{a,i}) + \eta^{i}, \label{forecast}
\end{eqnarray}
where $\eta^{i} \sim {\cal N}(0 , {\bf{Q}}_{t})$.
Upon receiving a new observation ${\bf{y}}_t$, each member of the forecast ensemble is adjusted using the Kalman gain ${\bf{K}}_{t}$ to generate the analysis ensemble ${\{ { {\bf{x}}}_{t}^{a,(i)} \}}_{i=1}^{N_{ens}}$ according to:

\begin{eqnarray}
{\bf{x}}^{a,i}_{t} &=& {\bf{x}}^{f,i}_{t} + {\bf{K}}_{t} ( {\bf{y}}^{i}_{t} -  \mathcal{H}_t {\bf{x}}^{f,i}_{t}), \\
{\bf{K}}_{t}  &=& {\bf{P}}_{t}^f \mathcal{H}_t^T  
(\mathcal{H}_t {\bf{P}}_{t}^f \mathcal{H}_t^T + {\bf{R}}_t)^{-1}, \label{K_gain}
\end{eqnarray}
where ${\bf{P}}_{t}^f$ denotes the sample forecast error covariance computed from the forecast members in \eqref{forecast} and ${\bf{y}}_{n}^{i}$ represents perturbed observations, i.e., ${\bf{y}}_{t}^{i} = {\bf{y}}_{t}+ \mu_{t}^{i}$ with  $\mu_{t}^{i}$ is a random noise sampled from the observational error distribution.

\section{Tracking Reference Solutions}
\label{results1}

\begin{figure}[!htbp]
    \centering
        
    \includegraphics[width=\linewidth]{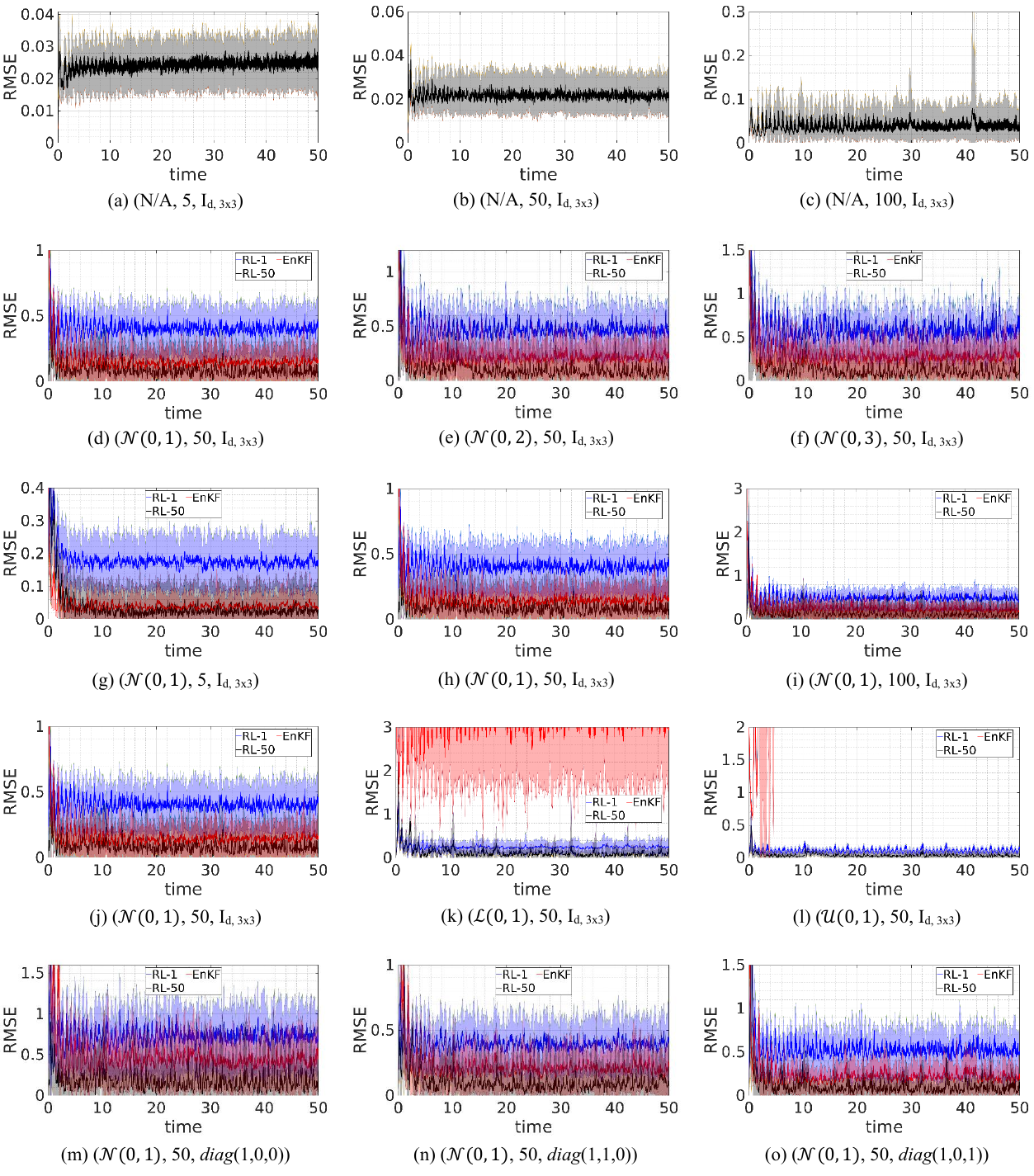}
        
    \caption{Evolution of the mean RMSE (solid lines) and its $\pm \sigma$ (shadowed) based on 50 experiment repetitions. Plotted are results for different experiments (a)-(c) tracking a noise-free reference solution, and for assimilating noisy observations in the case of (d)-(f) varying noise levels using normally-distributed noise, (e)-(i) different assimilation window lengths, (j)-(l) different noise distributions and (m)-(o) partial observability. The captions beneath each subplot describes the experimental condition in the order of noise distribution, $\delta t^o/\delta t$ the observation frequency and $\mathcal{H}$ the observation operator.}
    \label{fig:CompiledErrorsPlot}
    
\end{figure}

The RL-DA framework is systematically assessed under different experimental conditions. 
In the first scenario, an RL agent was trained to track a reference solution using coarse-in-time, noise-free observations of all state variables. 
Given the stochastic nature of the agent's policy function, the assimilated solution was not expected to precisely match the observations. 
Rather, the objective here was to investigate whether the corrections could maintain a reasonably close solution in comparison to the reference, and prevent them from diverging.
Three training regimes were explored, involving observations every 5, 50, and 100 $\delta t$, corresponding to $\delta t^o$ of 0.005, 0.05, and 0.1 time units, respectively. 
Evolution curves of the RMSEs of the RL solutions are presented in the top row of Figure \ref{fig:CompiledErrorsPlot}. 
The average RMSE is represented by a solid black line, encircled by a shaded region denoting one standard deviation ($\pm \sigma$), based on 50 repetitions of the experiment involving different reference solutions.
The plots indicate that the RMSE is on average approximately 0.025 for an assimilation frequency $\mathcal{T} = \delta t^o /\delta t$ values of 5 and 50, and increase to 0.05 for $\mathcal{T} = 100$. 
Furthermore, the top row of Figure \ref{fig:Compiled_z_Plot} illustrates RL and reference solutions for the $z$-variable in the Lorenz '63 system, based on randomly selected reference trajectories. 
These curves highlight strong agreement between the RL solution and the reference, further corroborating the results presented in Figure \ref{fig:CompiledErrorsPlot}.

\section{Assimilating Noisy Observations}
\label{results2}

In a more realistic scenario, an ensemble of noisy observations are assimilated to improve the model forecast. 
This investigation explores the influence of noise levels ($\sigma$), $\mathcal{T}$, statistical noise distribution, and partial state observability on the RL agent's performance.
Moreover, we conduct a comparative analysis by benchmarking the outcomes of the RL approach with those of the EnKF, which assimilates data from a relatively large ensemble comprising 50 realizations.
To ensure robustness and statistical significance, each of the RL and EnKF experiments was repeated 50 times using different reference solutions, providing a statistically significant estimate of the RMSE.

\subsection{Noise Level}

We examine the scenario of fully observed state available at regular intervals of $\mathcal{T} = 50$, with additive noise drawn from a Gaussian distribution characterized by zero mean and standard deviation $\sigma$. 
We investigate the influence of varying $\sigma$ on the assimilated solution by computing the RMSE for the complete trajectory, encompassing both forecast and analysis phases. 
We compare the results obtained from a single RL agent, an average solution derived from 50 distinct RL trajectories with actions randomly sampled from the agent's policy, and the EnKF solution based on an ensemble comprising 50 realizations.
Note that this comparison places the RL agent at a slight disadvantage, as it was trained without any statistical information about the response of the system to observation noise. 
Nonetheless, we believe that the comparison with the EnKF prediction is meaningful as it represents the primary benchmark against which DA algorithms are evaluated, despite the more suitable comparison with the Kalman Filter.
Notably, our algorithm consistently outperforms the Kalman Filter across all experiments and hence not shown.

The second row of Figure \ref{fig:CompiledErrorsPlot} presents the RMSE evolution over time for the assimilated solution, resulting from RL and EnKF under different $\sigma$ values. 
The plots suggest that, across all $\sigma$ values considered, the EnKF solution exhibits slightly lower RMSE values than those of a single RL agent, and slightly larger RMSEs than the RL solution obtained by averaging 50 action realizations.
This observation yields two significant insights: firstly, the potential computational efficiency gain from employing a single RL agent for DA, reducing computational overhead by a factor of at least $N_{ens}$, where $N_{ens}$ represents the ensemble size. 
Secondly, using a single RL agent with a stochastic policy allows for sampling a diverse set of forecast corrections, yielding a new ensemble of state estimates that when averaged, generally results in a lower RMSE compared to an EnKF solution produced using an equivalent ensemble size.

\begin{figure}[!htbp]
    \centering
    
    \includegraphics[width=\linewidth]{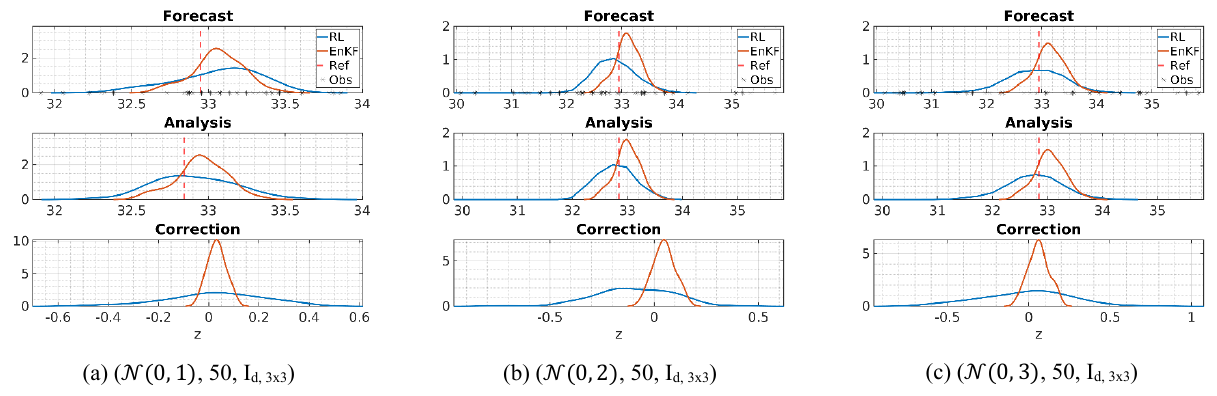}

    \caption{PDFs of the $z$-variable before (top) and after (middle) the correction step at time $t=45$ alongside the PDF of the correction (bottom) for the EnKF and RL solutions. The plots are presented for the experiment analyzing the sensitivity of the data assimilation algorithms to noise level.}
    \label{fig:PDFs_NoiseLevels}
    
\end{figure}

Figure \ref{fig:PDFs_NoiseLevels} illustrates the transition of the PDF after the correction is made alongside the distribution of the corrections for the RL and EnKF.
The results indicate that the RL distribution is wider and covers more of the observations points than the EnKF, meaning that the RL ensemble is richer in terms of information it provides even though individual realizations perform poorer than the EnKF solution. 
On the other hand, the mean of the RL solutions is closer to the reference solution than the average EnKF solution, aligning well with the results obtained earlier.
The plot also shows the distribution of the corrections, indicating that the distribution for the RL corrections is wider than that of the EnKF and suggesting that the EnKF is conservative when performing updates.
Similar results for the remaining experiments are analyzed in the Supplementary.

\begin{figure}[!htbp]
    \centering

    \includegraphics[width=\linewidth]{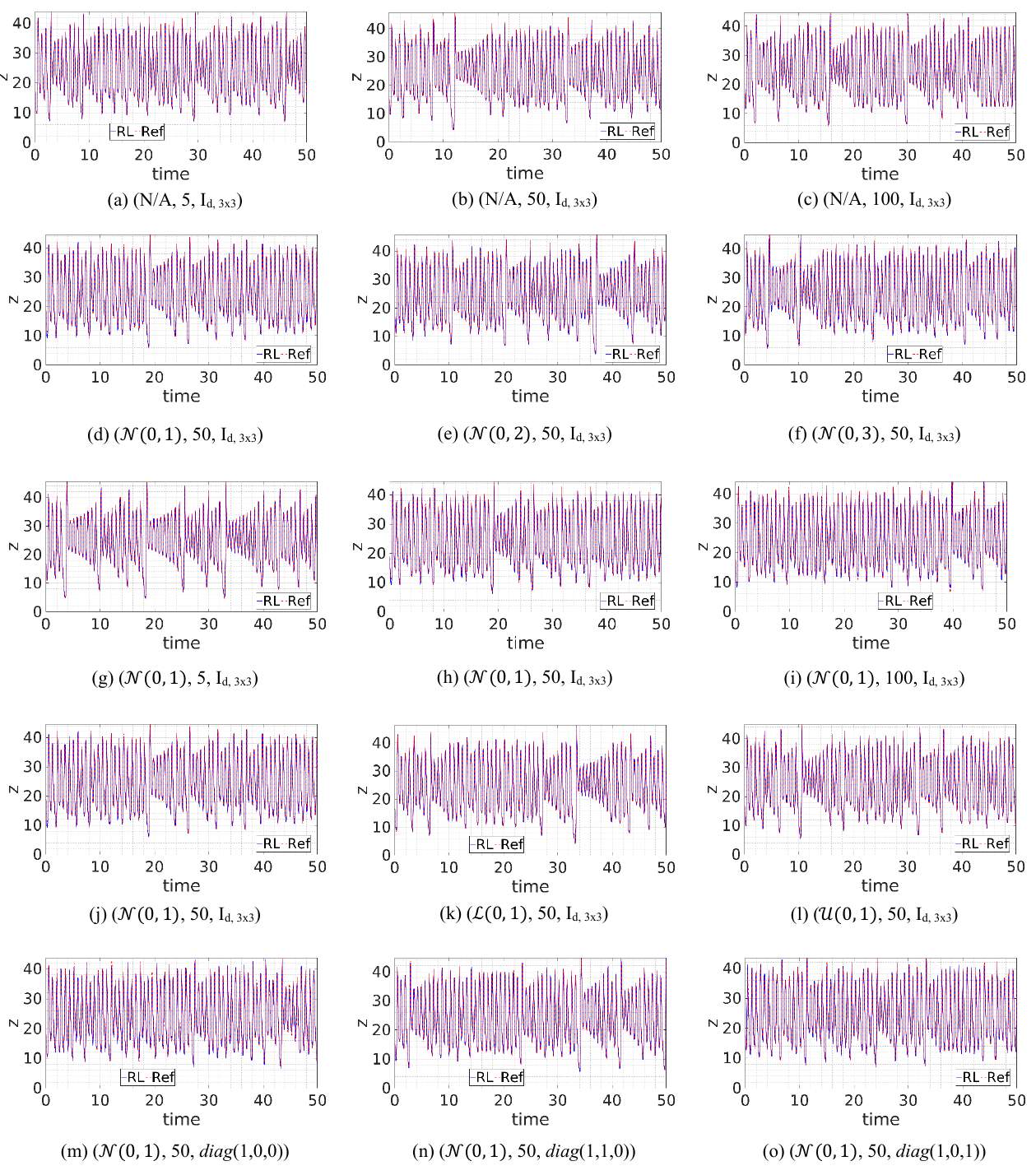}
        
    \caption{Evolution of the $z$-variable for a sample RL solution (solid blue lines) and corresponding reference (dashed red line). Plotted are results for different experiments (a)-(c) tracking a noise-free reference solution, and for assimilating noisy observations in the case of (d)-(f) varying noise levels using normally-distributed noise, (e)-(i) different assimilation window lengths, (j)-(l) different noise distributions and (m)-(o) partial observability. The captions beneath each subplot describes the experimental condition in the order of noise distribution, $\mathcal{T}$ the observation frequency and $\mathcal{H}$ the observation operator.}
    \label{fig:Compiled_z_Plot}
    
\end{figure}

As $\sigma$ increases, noticeable high-amplitude, abrupt variations in RMSE are observed in the assimilated solutions, and the time-averaged RMSE increases. 
In the second row of Figure \ref{fig:Compiled_z_Plot}, we present the RL and reference evolution curves corresponding to the $z$-variable. 
The results demonstrate that the RL solution closely follows the reference solution for all $\sigma$ values considered. 
However, as $\sigma$ increases, slight deviations between the RL solution and the reference become evident, particularly at the peaks and troughs of the curves. 
Nevertheless, the RL agent successfully assimilates noisy data, at high noise levels.

\subsection{Assimilation Frequency}

Observational data may often become available at varying time frequencies, necessitating a DA scheme capable of accommodating different observation rates. 
In light of this requirement, we trained an RL agent to assimilate noisy data for distinct $\mathcal{T}$, thereby examining the influence of high-frequency ($\mathcal{T} = 5$), medium-frequency ($\mathcal{T} = 50$), and low-frequency ($\mathcal{T} = 100$) observations.
The middle row of Figure \ref{fig:CompiledErrorsPlot} depicts the progression of RMSE under varying $\mathcal{T}$. 
Across all considered $\mathcal{T}$, the results suggest that a single RL agent exhibits slightly larger RMSE compared to those achieved by the 50-member EnKF solution.
For all cases, the 50 RL agent-averaged solution demonstrates a lower time-averaged RMSE in contrast to the 50-member averaged EnKF solution. 
This indicates that even when the RL agents do not communicate among each other, an MC averaged solution achieves lower RMSEs than the EnKF solution with 50 members.
Nevertheless, these results underscore the need to develop more sophisticated RL approaches, potentially utilizing multi-agent RL \cite{marl-book}, that incorporate ensemble information when performing the correction step.

\subsection{Noise Distribution}

A major limitation of the EnKF is its reliance on normally-distributed observations of system states. 
We investigate the impact of different statistical distributions of observations on the DA performance of the RL agent. 
Specifically, we examine cases involving unbiased standard Gaussian, strongly positively biased standard log-normal, and weakly positively biased standard uniform observational noise.
The $4^{th}$ row of Figure \ref{fig:CompiledErrorsPlot} presents the evolution curves of the RMSE for various observational noise distributions. 
The plots illustrate that for the case of standard Gaussian noise, both the single RL agent and EnKF solutions effectively assimilate noisy observational data with a slightly lower RMSE value achieved by the EnKF solution.
On the other hand, the 50-realization averaged RL solution yields a lower RMSE compared to the 50-member EnKF solution.
For log-normal and uniform noise distributions, the EnKF experiences large errors when assimilating noisy observations. 
Conversely, a single RL agent successfully assimilates these noisy observations, providing an assimilated solution that is close to the reference solution. 
Further improvements are observed when averaging the solutions obtained through policy sampling across 50 different realizations.
The penultimate row of Figure \ref{fig:Compiled_z_Plot} presents the RL and reference evolution curves for the $z$-variable. 
The plots indicate that the RL solution follows the reference solution reasonably well for all the noise distributions that were considered.
The curves clearly illustrate that the RL agent is able to assimilate non-Gaussian noisy observations even when observations are perturbed with biased noise.

\subsection{Partial Observability}

The practicality of DA lies in its ability to assimilate observations that partially or even indirectly characterize the evolution of state variables within a dynamical system. 
This is particularly valuable when the full system state cannot be directly observed, such as in real-world climate and weather applications.
To examine this setting, an RL agent was trained to assimilate noisy observations of select state variables--specifically, the $x$-variable alone, the $x$- and $y$-variables, and the $x$- and $z$-variables. 
The final row of Figure \ref{fig:CompiledErrorsPlot} portrays the evolution of RMSE of the aforementioned experiments.
The curves demonstrate that, in all cases, the RL agent provides a suitable correction that adequately guides the evolution of the full state. 
It is noteworthy that the RMSE of the solution obtained using a single RL agent is comparable to, albeit slightly higher than that of the EnKF with an ensemble of 50 realizations.
As observed in previous experiments, the averaged RL solution exhibits a lower average RMSE compared to the EnKF. 
To provide a tangible illustration of the assimilated solution's behavior, the final row of Figure \ref{fig:Compiled_z_Plot} presents curves depicting the temporal evolution of the $z$-variable for the case with partial system states observability. 
These plots depict that the RL assimilated solution generally tracks the reference, with occasional discrepancies that typically occur at the peaks and troughs, as expected.

\section{Discussion}
\label{conclusion}

This paper introduces RL as a novel approach for learning DA corrections. 
Through extensive experimentation on the Lorenz '63 dynamical system across various scenarios, we showcase the potential of the proposed approach. 
Our investigation encompasses both deterministic and stochastic settings, where RL agents are adeptly trained to track reference solutions and assimilate noisy data under varying conditions of assimilation window lengths, observational noise distributions, noise levels, and observed state variables.

The proposed RL-DA framework offers a paradigm shift by introducing new degrees of freedom to forecast-correction schemes, allowing for a nonlinear update term that satisfies a predefined optimal criteria, such as minimizing the root-mean-squared error in this study, hence, facilitating the discovery of novel correction strategies that are informed by the dynamical system through agent-environment interaction experiences.
Furthermore, RL imparts robustness to correction strategies, rendering them stable even in the presence of noisy perturbations and compounding errors. 
In this work, the RL agent minimizes the $\ell_2$ norm of the innovation term, a formulation demonstrated to be equivalent to maximizing the mutual information between observed state variables and their forecast counterparts.
Notably, this framework eliminates the need for a reference database as opposed to supervised learning approaches, which are commonly established through the assimilation of noisy observational data using methods such as the EnKF or variational methods \cite{Talagrand87}.

However, incorporating RL into DA raises critical questions warranting further exploration. 
While we employed the negative of the $\ell_2$ norm of the innovation term as the reward function in this study, more sophisticated functions considering system dynamics or ensemble information could potentially enhance the RL agent's performance. 
Moreover, since the RL agent is trained using the system of differential equations describing the evolution of a dynamical system, we speculate that this would force the agent to adapt and overcome model errors, when present.  
An overarching concern pertains to the physical validity of RL-derived solutions, which remains an open, fundamental question as is the case with other data-driven approaches when applied to physics-based applications. 
Although we did not directly encounter violations of physical constraints in our present setup, this avenue remains unexplored and in need for further exploration.

\section{Open Research}

All software and data used in the study will be made available upon acceptance at \url{https://github.com/mhammoud115/DA-RL}.


\acknowledgments
Research reported in this publication was supported by the Office of Sponsored Research (OSR) at King Abdullah University of Science and Technology (KAUST) CRG Award \#CRG2020-4336 and Virtual Red Sea Initiative Award \#REP/1/3268-01-01. The work of E.S.T. was supported in part by NPRP grant \# S-0207-200290 from the Qatar National Research Fund (a member of Qatar Foundation).

\clearpage
\newpage

\bibliography{agusample}

\clearpage
\newpage

\end{document}


\title{Data Assimilation in Chaotic Systems Using Deep Reinforcement Learning}

%
%

\authors{Mohamad Abed El Rahman Hammoud\affil{1}, Naila Raboudi\affil{1}, Edriss S. Titi\affil{2,3}, Omar Knio\affil{1} and Ibrahim Hoteit\affil{1}}

\affiliation{1}{King Abdullah University of Science and Technology, Thuwal 23955, Saudi Arabia}
\affiliation{2}{Department of Applied Mathematics and Theoretical Physics, University of Cambridge, Cambridge CB3 0WA, UK}
\affiliation{3}{Department of Mathematics, Texas A \& M University, College Station, TX 77843, USA}

\correspondingauthor{Ibrahim Hoteit}{ibrahim.hoteit@kaust.edu.sa}

\bigskip
\section*{This PDF file includes:}
\begin{list}{}{%
\setlength\leftmargin{2em}%
\setlength\itemsep{0pt}%
\setlength\parsep{0pt}}
\item Table S1
\item Figs.~S1 to S4
\end{list}

\newpage
\onecolumn

\begin{table}[!htbp]
\centering
\resizebox{\textwidth}{!}{%
\begin{tabular}{|l|l|l|l|l|}
\hline
                                                          & $\gamma$ & \text{max grad norm} & $v_f$ & $n_{a, train}$ \\ \hline
$\left( N/A,\ 5,\ I_{d, 3\times3} \right)$               & 0.9  & 0.9  & 0.7  & 100  \\ \hline
$\left( N/A,\ 50,\ I_{d, 3\times3} \right)$              & 0.1  & 0.8  & 0.7  & 100  \\ \hline
$\left( N/A,\ 100,\ I_{d, 3\times3} \right)$             & 0.1  & 0.9  & 0.7  & 100  \\ \hline
$\left( \mathcal{N}(0,1),\ 50,\ I_{d, 3\times3} \right)$ & 0.9  & 0.95 & 0.95 & 1000 \\ \hline
$\left( \mathcal{N}(0,2),\ 50,\ I_{d, 3\times3} \right)$ & 0.05 & 0.8  & 0.7  & 1000 \\ \hline
$\left( \mathcal{N}(0,3),\ 50,\ I_{d, 3\times3} \right)$ & 0.1  & 0.9  & 0.9  & 1000 \\ \hline
$\left( \mathcal{N}(0,1),\ 5,\ I_{d, 3\times3} \right)$  & 0.25 & 0.8  & 0.7  & 100  \\ \hline
$\left( \mathcal{N}(0,1),\ 50,\ I_{d, 3\times3} \right)$ & 0.9  & 0.95 & 0.95 & 1000 \\ \hline
$\left( \mathcal{N}(0,1),\ 100,\ I_{d, 3\times3} \right)$ & 0.05     & 0.95                  & 0.9                        & 1000                       \\ \hline
$\left( \mathcal{N}(0,1),\ 50,\ I_{d, 3\times3} \right)$ & 0.9  & 0.95 & 0.95 & 1000 \\ \hline
$\left( \mathcal{L}(0,1),\ 50,\ I_{d, 3\times3} \right)$ & 0.8  & 0.85 & 0.95 & 100  \\ \hline
$\left( \mathcal{U}(0,1),\ 50,\ I_{d, 3\times3} \right)$ & 0.1  & 0.9  & 0.8  & 100  \\ \hline
$\left( \mathcal{N}(0,1),\ 50,\  diag(1, 0, 0) \right)$  & 0.25 & 0.8  & 0.8  & 500  \\ \hline
$\left( \mathcal{N}(0,1),\ 50,\ diag(1, 1, 0) \right)$   & 0.3  & 0.9  & 0.7  & 500  \\ \hline
$\left( \mathcal{N}(0,1),\ 50,\ diag(1, 0, 1) \right)$   & 0.25 & 0.8  & 0.95 & 1000 \\ \hline
\end{tabular}%
}
\caption{Table describing the hyperparameters used to train the RL agent using the proximal policy optimization algorithm. The table outlines the hyperparameters for all 15 experiments considered in the study. All agents were trained using the ADAM optimization algorithm with a learning rate of $10^{-3}$. Moreover, all actor and critic networks are comprised of densely connected multilayer perceptrons with two hidden layers with 128 neurons each. }
\label{tab:my-table}
\end{table}


\clearpage
\newpage

\begin{figure}[!htbp]
    \centering
    
    \subfloat[$\left( \mathcal{N}(0,1),\ 5,\ I_{d, 3\times3} \right)$]{\includegraphics[width=0.32\linewidth]{figures/assimWindow/ev5/assimWindow_ev5_ev5_45.png}} \,
    \subfloat[$\left( \mathcal{N}(0,1),\ 50,\ I_{d, 3\times3} \right)$]{\includegraphics[width=0.32\linewidth]{figures/assimWindow/ev50/assimWindow_ev50_45.png}} \,
    \subfloat[$\left( \mathcal{N}(0,1),\ 100,\ I_{d, 3\times3} \right)$]{\includegraphics[width=0.32\linewidth]{figures/assimWindow/ev100/assimWindow_ev100_ev100_45.png}} \,

    \caption{PDFs of the $z$-variable before (top) and after (middle) the correction step at time $t=45$ alongside the PDF of the correction (bottom) for the EnKF and RL solutions. The plots are presented for the experiment analyzing the sensitivity of the data assimilation algorithms to assimilation frequency. As can be seen from the plots, the RL distributions for the $z$-variable and the correction are wider than that of the EnKF. Nevertheless, the RL distribution covers more of the noisy observations than the EnKF does. Furthermore, the mode of the RL ensemble is closer to the reference solution in comparison to the EnKF.}
    \label{fig:PDFs_assimWindow}
\end{figure}

\begin{figure}[!htbp]
    \centering
    \subfloat[$\left( \mathcal{N}(0,1),\ 50,\ I_{d, 3\times3} \right)$]{\includegraphics[width=0.32\linewidth]{figures/noiseDist/gaussian/Distribution_Gaussian_ev50_45.png}} \,
    \subfloat[$\left( \mathcal{L}(0,1),\ 50,\ I_{d, 3\times3} \right)$]{\includegraphics[width=0.32\linewidth]{figures/noiseDist/lognormal/Distribution_logNormal_ev50_45.png}} \,
    \subfloat[$\left( \mathcal{U}(0,1),\ 50,\ I_{d, 3\times3} \right)$]{\includegraphics[width=0.32\linewidth]{figures/noiseDist/uniform/Distribution_uniform_ev50_45.png}} \,
    \caption{PDFs of the $z$-variable before (top) and after (middle) the correction step at time $t=45$ alongside the PDF of the correction (bottom) for the EnKF and RL solutions. The plots are presented for the experiment analyzing the sensitivity of the data assimilation algorithms to the distribution of the observational noise. The plots indicate that while both the EnKF and RL distributions admit a high probability near the reference solution for the case of Gaussian noise, the RL solution is much closer to the reference solution than the EnKF solution in the case of lognormal and uniform noise. Furthermore, in the case of nongaussian noise, the EnKF correction term appears to be much more aggressive than that of RL and generally appears not to have a particular structure. }
    \label{fig:PDFs_distribution}
\end{figure}

\begin{figure}[!htbp]
    \centering
    
    \subfloat[$\left( \mathcal{N}(0,1),\ 50,\  diag(1, 0, 0) \right)$]{\includegraphics[width=0.32\linewidth]{figures/partial/X/partialObs_X_ev50_45.png}} \,
    \subfloat[$\left( \mathcal{N}(0,1),\ 50,\ diag(1, 1, 0) \right)$]{\includegraphics[width=0.32\linewidth]{figures/partial/XY/partialObs_XY_ev50_45.png}} \,
    \subfloat[$\left( \mathcal{N}(0,1),\ 50,\ diag(1, 0, 1) \right)$]{\includegraphics[width=0.32\linewidth]{figures/partial/XZ/partialObs_XZ_ev50_45.png}} \,

    \caption{PDFs of the $z$-variable before (top) and after (middle) the correction step at time $t=45$ alongside the PDF of the correction (bottom) for the EnKF and RL solutions. The plots are presented for the experiment analyzing the sensitivity of the data assimilation algorithms to partial observability. For the case of $\mathcal{H} = (1, 0, 0)$ and $(1, 1, 0)$, the plots indicate that the RL and EnKF distributions are comparable, where both cover most of the noisy observations and have the mode of the distribution close to the reference solution. Whereas for $\mathcal{H} = (1, 0, 1)$, the RL distribution is wider covering more of the noisy observations, and has the mode of the distribution closer to the reference solution in comparison to the EnKF solution.}
    \label{fig:PDFs_observations}
\end{figure}

\begin{figure}[!htbp]
    \centering
    \includegraphics[width = 0.9\linewidth]{figures/Sigma1_RMSE_vs_Ne.png}
    \caption{Plot illustrating the RMSE of the ensemble averaged solution as a function of the ensemble size $N_{ens}$. 
    The plot indicates that the RMSE of the EnKF solution saturates at an ensemble size of 10 meaning that an ensemble size of 50 is considered as a large cardinality ensemble for the Lorenz '63 system. On the other hand, the RMSE of the RL solution appears to keep on decreasing as $N_{ens}$ increases, with a much lower RMSE for small ensembles. This suggests that the RL framework offers huge computational savings with an adequately reliable solution, especially when computational resources are scarse.}
    \label{fig:ensSize}
\end{figure}